\documentclass[sigconf]{acmart}
\usepackage[utf8]{inputenc}
\usepackage{libertine}
\usepackage{algorithmic}
\usepackage{graphicx}
\usepackage{textcomp}
\usepackage{xcolor}
\usepackage{stmaryrd}
\usepackage{caption}
\usepackage{subcaption}
\usepackage{siunitx}
\usepackage{anyfontsize}

\AtBeginDocument{%
  }

\acmConference[PASC '25]{}{June 16--18,
  2025}{Brugg, Switzerland}

\newcommand{\rev}[1]{#1}

\begin{document}
\title[Space-time parallel scaling of Parareal with a PINO coarse propagator for the Black-Scholes equation]{Space-time parallel scaling of Parareal with a \rev{physics-informed} Fourier Neural Operator coarse propagator \rev{applied to the Black-Scholes equation}}

\author{Abdul Qadir Ibrahim}
\email{abdul.ibrahim@tuhh.de}
\orcid{0000-0003-3452-8500}
\affiliation{%
  \institution{Chair Computational Mathematics, Institute of Mathematics}
  \city{Hamburg}
  \country{Germany}
}

\author{Sebastian G\"otschel}
\email{sebastian.goetschel@tuhh.de}
\orcid{0000-0003-0287-2120}
\affiliation{%
  \institution{Chair Computational Mathematics, Institute of Mathematics}
  \city{Hamburg}
  \country{Germany}}

\author{Daniel Ruprecht}
\email{daniel.ruprecht@tuhh.de}
\orcid{0000-0003-1904-2473}
\affiliation{%
  \institution{Chair Computational Mathematics, Institute of Mathematics}
  \city{Hamburg}
  \country{Germany}
}

\renewcommand{\shortauthors}{Ibrahim et al.}
\begin{abstract}
  Iterative parallel-in-time algorithms like Parareal can extend scaling beyond the saturation of purely spatial parallelization when solving initial value problems.
    However, they require the user to build coarse models to handle the unavoidable serial transport of information in time.
    This is a time-consuming and difficult process since there is still limited theoretical insight into what constitutes a good and efficient coarse model.
    Novel approaches from machine learning to solve differential equations could provide a more generic way to find coarse-level models for multi-level parallel-in-time algorithms.
    This paper demonstrates that a physics-informed Fourier Neural Operator (PINO) is an effective coarse model for the parallelization in time of the two-asset Black-Scholes equation using Parareal.
    We demonstrate that PINO-Parareal converges as fast as a bespoke numerical coarse model and that, in combination with spatial parallelization by domain decomposition, it provides better overall speedup than both purely spatial parallelization and space-time parallelization with a numerical coarse propagator.
\end{abstract}



\maketitle

\section{Introduction}\label{sec:introduction}
Space-time parallelization, combining spatial parallelization via domain decomposition with parallel-in-time integration, has been shown to be able to extend parallel scaling beyond what parallelization in space alone can provide~\cite{SpeckEtAl2012}.
However, a key challenge when deploying iterative, multi-level methods like Parareal~\cite{LionsEtAl2001}, MGRIT~\cite{FalgoutEtAl2014_MGRIT} or PFASST~\cite{EmmettMinion2012} remains the need to define one or multiple coarse level models that deal with the inevitable serial information propagation in the time direction.
A good coarse-level model must be accurate enough to ensure quick convergence of the method but also needs to be computationally cheap, as it forms a serial bottleneck and limits achievable speedup.
Physics-informed neural operators~\cite{LiEtAl2021_pino}, an approach from machine learning to solve partial differential equations, have two properties that suggest they might be attractive as coarse-level models for parallel-in-time methods: first, they are very fast to evaluate once trained and second they are generic, requiring just the PDE residual in the loss function, although some tuning of training parameters might be needed.
Neural operators can solve initial-boundary value problems over some time horizon in a single step. 
Together with an appropriate discretization of the PDE in space, the Neural Operator can be trained to map the state at time $t$ to the state at a later point $t + \Delta t$~\cite{GradyEtAl2023}. 
Hence, the neural operator performs exactly the mapping required from a coarse model in Parareal.
In contrast to sparse, mesh-based numerical algorithms, neural networks also run efficiently on GPUs, suggesting that a combination of both could help to better utilise heterogeneous architectures~\cite{IbrahimEtAl2023}.
We investigate the effectiveness of a physics-informed neural operator (PINO) as coarse level model for the parallel-in-time integration method Parareal to solve the two-asset Black-Scholes equation, a PDE used in computational finance.
The key novelty in this paper is the demonstration that Parareal with an ML-based coarse model not only provides speedup over serial time stepping but can extend scaling beyond the saturation point of space-only parallelization \rev{and that a PINO provides a more effective coarse model than a PINN-based propagator.}

\section{Related Work}\label{sec:related-work}
This paper is an extension of Ibrahim et al. 2023~\cite{IbrahimEtAl2023} where we study a physics-informed neural network propagator (PINN-P) for the one-asset Black-Scholes equation as a coarse model for Parareal. 
The PINN-P propagator was inspired by classical physics-informed neural networks (PINNs), but, instead of solving the differential equation directly, was used to learn the mapping from a given initial value to the solution of the Black-Scholes equation at a future time step.
PINN-P therefore used a PINN-like architecture to create a simple Neural Operator (NO). 
NOs, like the Fourier Neural Operator (FNO)~\cite{LiEtAl2021_fno} or DeepONets~\cite{LuEtAl2021}, are methods to learn such solution operators, i.e., the mapping from parameters or initial conditions to the solution of a partial differential equation. 
This is in contrast to an actual PINN that would learn the solution of a PDE for a fixed initial condition. 

Here, we extend our previous results in three ways: (i) we show that a Physics Informed Neural Operator (PINO), based on the Fourier neural operator, is a better coarse propagator than PINN-P as it requires fewer parameters and a significantly shorter training time, 
(ii) we consider the more complex two-asset Black-Scholes equations and combine Parareal with spatial parallelization, (iii) we analyze performance of the full space-time parallelization in contrast to studying the performance of Parareal alone.

We have discussed the main works studying combinations of machine learning with parallel-in-time integration in more detail in our previous paper~\cite{IbrahimEtAl2023}.
Therefore, we only provide a quick summary whilst discussing publications that have been published since then in more detail below.
Using ML to build coarse propagators was first studied by Yalla and Enquist~\cite{YallaEnquist2018}.
The first reports of speedups were made by Agboh et al.~\cite{AgbohEtAl2020} for a system of ordinary differential equations related to robotic control.
Nguyen and Tsai use ML not to replace the coarse model but to enhance it, in order to improve its accuracy~\cite{NguyenEtAl2023b}.
Gorynina et al. use ML to provide fast approximations of force fields in molecular dynamics simulations with Parareal~\cite{GoryninaEtAl2022}. 

A few relevant papers on the combination of data-based techniques with Parareal have been published since our previous paper.
Motivated by machine learning approaches, Jin et al.~\cite{JinEtAl2023} recently proposed an optimization method to construct optimized coarse propagators, based on error estimates, that exhibit favourable numerical properties.
Focused on wave-propagation problems, Kaiser et al.~\cite{kaiser2024efficient} extend their previous work~\cite{NguyenEtAl2023b} to provide a more general end-to-end framework. 
Instead of learning a coarse propagator, they combine a coarse numerical solver with a convolutional auto-encoder (JNet) to augment the under-resolved wave fields.
Their approach tries to improve coarse level accuracy to get Parareal to converge in fewer iterations but not to reduce the computational cost of the coarse level.
By contrast, we aim to reduce the computational cost of the coarse model to improve the scaling of Parareal but not necessarily to reduce the number of iterations.
An approach similar to the one by Kaiser et al. is pursued by Fang and Tsai~\cite{FangEtAl2023b} to stabilize Parareal for Hamiltonian systems.
They suggest using either a neural network to correct phase information in the coarse model or to replace the fine model.
Both approaches are tested for the Fermi-Pasta-Ulam problem. A combination of Parareal with a Neural Operator for fusion MHD simulations was very recently investigated~\cite{PamelaEtAl2025} but without analyzing speedups from space-time parallelization.
Betcke et al. use a Random Projection Neural Network as coarse propagator for Parareal~\cite{BetckeEtAl2024}.
They provide a theoretical analysis of convergence properties of the resulting method but do not measure speedups.
Gattiglio et al.~\cite{GattiglioEtAl2024b} propose a combination of Random Neural Networks with Parareal called RandNet-Parareal.
They show a significant improvement in speedup over classical Parareal, in line with PINN-Parareal~\cite{IbrahimEtAl2023}, but do not explore space-time parallelization.

While our approach leverages machine learning to improve parallel-in-time integration algorithms, there are also studies that go the other way, using Parareal to improve training of neural networks.
Lee et al.~\cite{LeeEtAl2022b} further develop their approach to use a Parareal-like technique to accelerate the training of very deep networks, similar to previous works by Guenther et al.~\cite{GuentherEtAl2020} or Meng et al.~\cite{MengEtAl2020}.
They interpret a deep neural network as what they call a ``parareal neural network'', which consists of coarse and fine structures, similar to the propagators in Parareal.
A Parareal-style iteration is then used to parallelize the training of the DNN across multiple GPUs. 

One issue when using machine learning to construct neural operators is the required amount of training data. 
To reduce this, physics-informed variants add the PDE residual and other constraints to the loss function, and thus can even be trained without any training data. 
Physics-informed neural operators~\cite{LiEtAl2021_pino} based on the Fourier neural operator~\cite{LiEtAl2021_fno}, as used here, were originally constructed to approximate PDE solution operators on cartesian, periodic domains, but can be extended to more general geometries~\cite{LiEtAl2023b}. 
Deep operator networks, so-called DeepONets~\cite{LuEtAl2021}, have been extended to physics-informed variants as well~\cite{GoswamiEtAl2023,WangEtAl2021}.  
Alternative architectures for operator learning include, among others, nonlinear manifold decoders~\cite{FangEtAl2024,SeidmanEtAl2022} and wavelet-based approaches~\cite{GuptaEtAl2021,NavaneethEtAl2024}.

We use the Black-Scholes equation from computational finance with two assets as a test problem~\cite{Merton1973,BlackScholes1973}. 
The performance of classical Parareal for option pricing has been studied before~\cite{BalMaday2002,PagesEtAl2018,GuEtAl2024} but not with ML-based coarse propagators.

\section{Algorithms and Benchmark Problem}\label{sec:algorithms-and-benchmark-problem}
Section~\S\ref{subsec:parareal} briefly explains the Parareal algorithm.
In \S\ref{subsec:mesh-based-discretization} we introduce our benchmark problem, the two-asset Black-Scholes equation, and describe a numerical solution approach based on finite differences.
Then, \S\ref{subsec:space-parallelization} describes the spatial parallelization while \S\ref{subsec:pino} explains the physics-informed neural operator we use to construct a coarse-level model for Parareal.
\subsection{Parareal}\label{subsec:parareal}
Parareal is an algorithm for the parallel solution of an initial value problem
\begin{equation}
    \label{eq:ivp}
u'(t) = f(u(t)), \ u(t_0) = u_0
\end{equation}
for $0 \leq t \leq T$ which here stems from the spatial discretization of a time-dependent PDE with finite differences, see \S\ref{subsec:mesh-based-discretization}.
To parallelize integration of~\eqref{eq:ivp},    Parareal needs a decomposition of the temporal domain $[0,T]$ into $P_{\text{time}}$ time slices $[T_n, T_{n+1}]$, $n=0,\ldots,P_{\text{time}}-1$.
Although it is possible to assign multiple time slices to a processor, we always assume that every processor handles only one-time slice so that $P_{\text{time}}$ is also equal to the number of processors used in the temporal direction.\footnote{Note that when Parareal is combined with spatial parallelization, every time slice is assigned not to a single processor but to $P_{\text{space}}$ many processors, where $P_{\text{space}}$ is the number of processors over which a single instance of the solution $u_n$ is distributed in space. The total number of processors is then $P_{\text{time}} \times P_{\text{space}}$.}
Let $\mathcal{F}$ denote a time integration method of suitable accuracy such that
\begin{equation}
    \label{eq:fine}
	u_{n+1} = \mathcal{F}(u_n)
\end{equation}
is the approximation delivered at the end of a slice $[T_n, T_{n+1}]$ when starting from an initial value $u_n$ at $T_n$.
Classical time stepping would compute~\eqref{eq:fine} sequentially for $n=0, \ldots, P_{\text{time}}-1$.
Parareal replaces this serial procedure with the iteration
\begin{equation}
\begin{aligned}
u^{k+1}_{n+1} &= \mathcal{G}(u^{k+1}_n) + \mathcal{F}(u^k_n) - \mathcal{G}(u^k_n), \\
    &\quad n=0, \ldots, P_{\text{time}}-1, \\
    &\quad k=0, \ldots, K-1
\end{aligned}
\label{eq:parareal}
\end{equation}
where $K$ is the number of iterations and $\mathcal{G}$ a coarse solver for~\eqref{eq:ivp}.
\rev{Note that both fine and coarse propagator in Parareal are treated as black-boxes and can encompass explicit or implicit methods as well as arbitrary direct or iterative linear or nonlinear solvers~\cite{Aubanel2011}.}

Although the serial computation of iteration~\eqref{eq:parareal} is computationally more expensive than serial time stepping by computing~\eqref{eq:fine} for $n=0, \ldots, P_{\text{time}-1}$, the costly evaluation of $\mathcal{F}$ can be parallelized over all time slices. 
Therefore, if the coarse model is cheap enough and the number of iterations small, iteration~\eqref{eq:parareal} can solve~\eqref{eq:ivp} in less wallclock time than serially evaluating~\eqref{eq:fine}.
For Parareal, speedup using $P_{\text{time}}$ timeslices is bounded by 
\begin{equation}
    \label{eq:speedup}
S(P_{\text{time}}) \leq \frac{1}{ \left(1 + \frac{K}{P_{\text{time}}} \right)  \frac{c_{\text{coarse}}}{c_{\text{fine}}} + \frac{K}{P_{\text{time}}}} \leq \min\left\{ \frac{c_{\text{fine}}}{c_{\text{coarse}}}, \frac{P_{\text{time}}}{K} \right\}
\end{equation}
where $c_{\text{coarse}}$ and $c_{\text{fine}}$ are the execution times for $\mathcal{G}$ and $\mathcal{F}$.
To allow for speedup from Parareal, $K$ has to be small, that is, the method needs to converge quickly and $c_{\text{coarse}} \ll c_{\text{fine}}$, that is, the coarse propagator must be computationally much cheaper than the fine.
We will demonstrate that a PINO provides a coarse propagator that is much faster than a numerical $\mathcal{G}$, thus reducing $c_{\text{coarse}}$, but accurate enough to not require more iterations $K$. 
It therefore provides a significant benefit in scaling in situations where the first term in~\eqref{eq:speedup} is limiting speedup.
Details on Parareal, including the pseudocode of a Parareal implementation, can be found in the literature~\cite{Ruprecht2017_lncs}.

\subsection{Numerical solution of the two-asset Black-Scholes equation}\label{subsec:mesh-based-discretization}
We consider the Black-Scholes equation in the form
\begin{align}
\label{eq:bse}
\frac{\partial u}{\partial t} 
&+ \frac{1}{2}\sigma_1^2 x^2 \frac{\partial^2 u}{\partial x^2} 
+ \frac{1}{2}\sigma_2^2 y^2 \frac{\partial^2 u}{\partial y^2} \notag \\
&+ \rho \sigma_1 \sigma_2 xy \frac{\partial^2 u}{\partial x \partial y} + rx \frac{\partial u}{\partial x} + ry \frac{\partial u}{\partial y} - ru = 0
\end{align}
where $u(t, x(t), y(t))$ is the option price, $t$ is time, $x$ the first and $y$ the second asset value, $\sigma_1$ and $\sigma_2$ are the volatilities of the two assets, $\rho$ is the correlation of the two assets and $r$ is the risk-free interest rate.
\rev{Note that while~\eqref{eq:bse} is a model from computational finance, its mathematical structure is that of a quasi-linear advection-diffusion-reaction equation. 
Therefore, the presented approach will hopefully generalize to other applications.}

Since the PDE is typically augmented with a final time condition, we change the variable $\tau  = T -t$ to turn that into an initial condition for convenience. 
Under this change of variable, the equation becomes
\begin{align}
\label{eq:bse_change_variable}
\frac{\partial u}{\partial \tau} 
&- \frac{1}{2}\sigma_1^2 x^2 \frac{\partial^2 u}{\partial x^2} 
- \frac{1}{2}\sigma_2^2 y^2 \frac{\partial^2 u}{\partial y^2} \notag \\
&- \rho \sigma_1 \sigma_2 xy \frac{\partial^2 u}{\partial x \partial y} 
- rx \frac{\partial u}{\partial x} 
- ry \frac{\partial u}{\partial y} 
+ ru = 0
\end{align}
This can be extended to a more general n-dimensional PDE
\begin{align}
	\label{eq:bse_general}
	\frac{\partial u}{\partial t} + \frac{1}{2}\sum_{i,j=1}^{n} \rho_{ij}\sigma_i\sigma_j x_i x_j \frac{\partial^2 u}{\partial x_i \partial x_j} + r \sum_{i=1}^{n} x_i \frac{\partial u}{\partial x_i} - r u = 0
\end{align}
for $n$ assets but we focus on the two-asset case with $n=2$ here.
We use a centered finite difference for spatial discretization and an implicit Euler scheme for time discretization~\cite{Yongsik2013}, resulting in the discretization
\begin{subequations}
\label{eq:fd}
\begin{align}
    &\frac{u_{i,j}^{n+1} - u_{i,j}^n}{\Delta \tau} 
    - \frac{1}{2}\sigma^{2}_1 x^2_{i} \frac{u_{i-1,j}^{n+1} - 2u_{i,j}^{n+1} + u_{i+1,j}^{n+1}}{\Delta x^2} \notag \\
    &- \frac{1}{2}\sigma^{2}_2 y^2_{j} \frac{u_{i,j-1}^{n+1} - 2u_{i,j}^{n+1} + u_{i,j+1}^{n+1}}{\Delta y^2} \notag \\
    &- \rho \sigma_1 \sigma_2 x_i y_j \frac{u_{i+1,j+1}^{n+1} + u_{i-1,j-1}^{n+1} - u_{i+1,j-1}^{n+1} - u_{i-1,j+1}^{n+1}}{4 \Delta x \Delta y} \notag \\
    &- r x_i \frac{u_{i+1,j}^{n+1} - u_{i-1,j}^{n+1}}{2 \Delta x} 
    - r y_j \frac{u_{i,j+1}^{n+1} - u_{i,j-1}^{n+1}}{2 \Delta y} 
    + r u_{i,j}^{n+1} = 0
\end{align}
\end{subequations}
The subscripts $i$ and $j$ stand for horizontal and vertical node index and the superscript $n$ for the time step index. 
For simplicity, we consider a uniform grid with $\Delta x = \Delta y \rev{ =1}$ where the nodal distances for the two axes are same. 
The time step size for the coarse propagator is chosen as \( \Delta \tau = 0.1 \). For the fine propagator, this interval is subdivided into three smaller steps.
We consider the two-asset cash-or-nothing opion since the closed form solution for this is known~\cite{math8030391} and can be used to verify the accuracy of PINO-Parareal.

\rev{The closed-form solution for the cash-or-nothing option reads:}
\begin{align}
	\label{eq:solution}
	& u(x,y,\tau) = ce^{-r\tau} B(d_x,d_y;\rho), \\
	&\quad d_x = \frac{\ln\left(\frac{x}{S_1}\right)+(r-0.5\sigma_x^2)\tau}{\sigma_x\sqrt{\tau}}, \\ 
	&\quad d_y = \frac{\ln\left(\frac{y}{S_2}\right)+(r-0.5\sigma_y^2)\tau}{\sigma_y\sqrt{\tau}},
\end{align}
where the Bivariate cumulative normal distribution function~\cite{Genz2004} is given by
\begin{align}
	B(d_x,d_y;\rho) = \frac{1}{2\pi\sqrt{1-\rho^2}}\int_{-\infty}^{d_x}\int_{-\infty}^{d_y}e^{-\frac{\xi_1^2-2\rho\xi_1\xi_2+\xi_2^2}{2(1-\rho^2)}}d\xi_2d\xi_1.
\end{align}

\rev{The computational domain is $[0, 300]\times[0, 300]$ for $t \in [0, 1]$ and~\eqref{eq:solution} is used to set the Dirichlet boundary condition.
The exercise price is set to $S_1 = 100$ and $S_2 = 100$, the risk-free interest rate to $r = 1.0$, the volatility of the first asset to $\sigma_1 = 0.3$,
the volatility of the second asset to $\sigma_2 = 0.3$, the correlation to $\rho = 0.5$ and the maturity is $T = 1$~\cite{math8030391}.
These parameters are chosen to balance practical relevance and the need to test Parareal in a regime where achieving good convergence is a challenge.
Therefore, the risk-free interest rate \( r = 1.0 \) chosen higher than typical market rates to magnify the drift terms in~\eqref{eq:bse_change_variable}, since these make achieving good Parareal convergence more difficult~\cite{GanderVandewalle2007_SISC}. 
Volatilities \( \sigma_1 = \sigma_2 = 0.3 \) are moderate and ensure that the diffusion terms contribute noticeably to the dynamics without making the problem diffusion-dominated, a regime in which Parareal is already known to converge very quickly.
A correlation of \( \rho = 0.5 \) introduces interdependence between the assets, which is crucial for capturing realistic multi-asset interactions. 
In real-world applications, these parameters depend on market conditions: volatilities typically range from 0.1 to 0.5, interest rates normally vary between 0.01 and 0.05, and correlations must be in \([-1, 1]\). 
The selected values test Parareal as well as the PINO model's robustness without being too unrealistic. 
}

\rev{Strike prices \( S_1 = S_2 = 100 \) represent an at-the-money option, which produces meaningful dynamics around the most sensitive region of the option price. 
The strike prices $S_1, S_2$ do not explicitly appear in the Black-Scholes equation but are incorporated into the terminal or boundary conditions.
For this problem, the terminal payoff function of cash-or-nothing options with wo assets is given by
\[
u_T(x, y) =
\begin{cases} 
c, & \text{if } x \geq S_1, \ y \geq S_2,  \\ 
0, & \text{otherwise.}
\end{cases}
\]
where $S_1$, $S_2$ are the strike prices of the two assets and $c$ is the cash amount received at expiration.
\begin{figure*}[!t]
\centering
\includegraphics{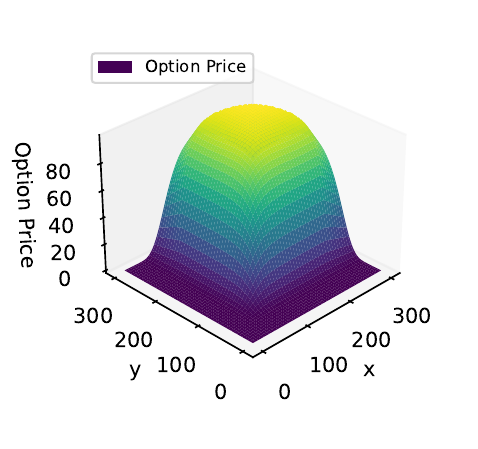}
\hfil
\includegraphics{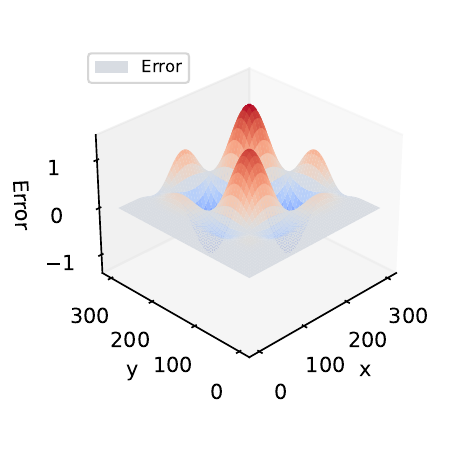}
\rev{\caption{Numerical solution using finite differences (left) of the two-asset Black-Scholes equation and error against the closed-form solution (right).}}
\end{figure*}
}

\subsection{Spatial parallelization}\label{subsec:space-parallelization}
We use the \texttt{dask} Python framework~\cite{dask} to parallelize the computation of the finite difference approximation~\eqref{eq:fd} in space.
This library offers high-level abstractions for distributed computing and uses \text{mpi4py} for communication and data exchange.
\rev{Note that while a shared memory parallelization in time might offer advantages, the global interpreter lock in Python makes this very difficult to realize.
We therefore rely on distributed memory parallelization in both space and time.}
Spatial parallelization is done by distributing the iterative solver step in each time step of implicit Euler by decomposing the spatial domain.
To do this, we implemented a conjugate gradient (CG) method using \texttt{dask} arrays to parallelize the matrix-vector and vector-vector products arising in CG, which are by far the most computationally costly part of the algorithm.
For spatial parallelization using $P_{\text{space}}$ processors, speedup is bounded by
\begin{equation}
 S(P_{\text{space}}) \leq P_{\text{space}}.
\end{equation}
For a combined space-time parallelization, in the best case the achieved speedup is multiplicative so that
\begin{equation}
 S(P_{\text{space}}, P_{\text{time}}) \leq  \frac{P_{\text{space}}}{ \left(1 + \frac{K}{P_{\text{time}}} \right)  \frac{c_{\text{coarse}}}{c_{\text{fine}}} + \frac{K}{P_{\text{time}}}}
\end{equation}
Note that the total number of processors is $P_{\text{total}} = P_{\text{space}} \times P_{\text{time}}$ since computations on every time slice are parallelized using $P_{\text{space}}$ many \rev{processes}.
Also, since this model neglects overheads from communication or competition for resources, it is only an upper bound and not necessarily predictive of actual speedup.

\rev{\subsection{Fourier Neural Operator (FNO)}\label{subsec:fno}
The Fourier Neural Operator (FNO)~\cite{LiEtAl2021_fno} is a recently developed deep learning framework for 
approximating solution operators of partial differential equations (PDEs). Unlike traditional numerical 
methods that compute the solution for a fixed discretization, FNO learns a mapping between function spaces, 
enabling generalization across different resolutions and input conditions. At its core, FNO leverages the Fast 
Fourier Transform (FFT) to express the integral kernel of the operator in the frequency domain. In each Fourier 
layer, the input function is first transformed to Fourier space, where a learnable multiplier is applied to a 
truncated set of Fourier modes. An inverse FFT then returns the modified data to physical space. 
This procedure allows the network to capture long-range dependencies and global features while keeping 
the number of learnable parameters modest. The typical FNO architecture involves three main steps:
\begin{itemize} 
	\item Lifting: The input function is projected into a higher-\-dimen\-sional feature space using a pointwise fully connected layer.
	\item Fourier Layers: A sequence of layers is applied where each layer computes the FFT of the feature maps, performs a modal truncation and a corresponding multiplication with trainable weights in Fourier space, and finally transforms the result back with an inverse FFT.
	\item Projection: The resulting representation is mapped back to the target function space via another pointwise layer.
\end{itemize}
Key advantages of the FNO framework include its mesh invariance, which permits the application of the trained 
operator to inputs with different discretizations, and its computational 
efficiency stemming from the FFT operations. In our work, the FNO serves as the backbone for 
the physics-informed neural operator (PINO). By embedding the physics of the Black-Scholes PDE into the loss 
function, the PINO retains the efficiency of FNO while ensuring that the learned operator respects the governing equations.}

\subsection{Physics Informed Neural Operator (PINO)}\label{subsec:pino}
We present a specific application of physics-informed neural operators for solving initial value problems (IVPs) using a Fourier Neural Operator (FNO) \cite{ZongyiEtAl202}. 
During model training, we randomly generate 10,000 collocation points uniformly within the spatial domain, 5,000 collocation points at the boundary and 5,000 collocation points at expiration.  
We also generate 5,000 initial conditions for the IVP. 
These initial conditions were sampled at different time steps of the equation to cover the range of possible states within the system's dynamics.
The loss function to be minimized is similar to the single-asset loss functions~\cite{IbrahimEtAl2023},
\begin{equation}
\text{MSE}_{\text{total}} = \text{MSE}_f + \text{MSE}_{\exp} + \text{MSE}_b,
\label{eq:bs_totalloss}
\end{equation}
consisting of a term to minimize the PDE residual $f(u)$
\begin{align}
\text{MSE}_f &=  \frac{1}{N_f}\sum_{i=1}^{N_f} |f(\tilde{u}(t_i, x_i, y_i))|^2,
\label{eq:loss_pde_call}
\end{align}
the boundary loss term
\begin{equation}
\text{MSE}_{\text{b}} = \frac{1}{N_b}\sum_{i=1}^{N_b}\left|\tilde{u}(t_i, x_i, y_i) - u(t_i, x_i, y_i)\right|^2,
\label{eq:loss_bc_call}
\end{equation}
and the terminal (expiration) loss 
\begin{equation}
\text{MSE}_{\exp} = \frac{1}{N_{\exp}}\sum_{i=1}^{N_{\exp}}\left|\tilde{u}(T, x_i, y_i) - \max[\max(x_i, y_i) - S. 0)]\right|^2,
\label{eq:loss_fs_call}
\end{equation}

This loss function is an approximation of the $L^2$-loss from Li et al.~\cite{LiEtAl2021_pino}.
The model is configured with several input arguments: the time interval for training, generated initial condition functions, the number of domain and boundary points, and the number of batches. 
The default configuration includes a Fourier neural operator backbone with a width (the number of channels in each hidden layer) of 64, 
a mode (the number of Fourier modes used) of 12, and L (the number of stacked Fourier integral operator layers) set to 4. 
Additionally, ReLU activation functions and batch normalization are applied. For optimization, we employ the Adam optimizer~\cite{kingma2017adammethodstochasticoptimization} to train the model for 2,500 epochs with an initial learning rate of 0.001, decay steps of 25 epochs and a decay rate of 0.96.
\rev{While using L-BFGS can provide benefits to training~\cite{ZampiniEtAl2024}, experiments not documented here show that, in our case, it only leads to slightly faster initial convergence but eventually delivers comparable accuracy as Adam at 2500 epochs.
Therefore, Adam leads to overall shorter training times here.}

\begin{figure*}[!t]
	\centering
	\includegraphics{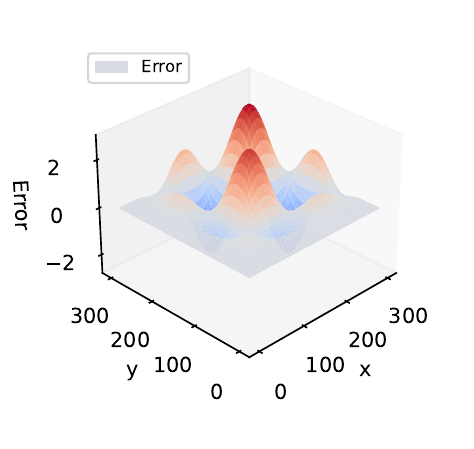}
	\hfil
	\includegraphics{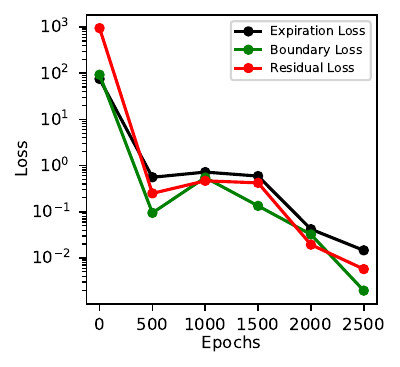}
	\caption{Error of the PINO-generated solution against the closed-form solution (left) and change of the PINO loss function during training (right). After 2500 epochs the loss function is reduced by five orders of magnitude.}
	\label{fig:combined_1}
\end{figure*}

To compute its accuracy, we compare the PINO solution to the closed-form solution~\eqref{eq:solution}.
To use the PINO as a coarse propagator in Parareal, the trained model takes as input the asset values at time $t_n$, the time interval ($t_n, t_{n+1}$), 
and returns the solution at the end point of the time interval $t_{n+1}$.  
The training time for the PINO is approximately 10 minutes.
By contrast, the PINN-P we used in earlier work~\cite{IbrahimEtAl2023}, modified to solve the two-asset Black-Scholes equation~\eqref{eq:bse_change_variable}, took more than 40 minutes to train.
Figure~\ref{fig:combined_1} (left) shows the solution produced by the PINO while Figure~\ref{fig:combined_1} (right) shows the three components of the loss function over the training process.

\section{Numerical Results}\label{sec:results}
The system used for the tests is the Linux cluster at Hamburg University of Technology (TUHH), running AlmaLinux 8 and equipped with fast OmniPath connectivity. 
The system architecture supports both 32-bit and 64-bit operation modes. 
It features different types of nodes but we use one with two AMD EPYC 9124 32-core CPUs. The CPU frequency ranges from 1500 MHz to 3711.9141 MHz, with a base 
frequency of 3000 MHz. 
The system node has 32 KB each of L1 data and instruction cache, 1024 KB of L2 cache, and 16384 KB of L3 cache. 
Additionally, the system is configured 
with NUMA (Non-Uniform Memory Access) architecture, consisting of two NUMA nodes. 
Virtualization support is enabled, and the system's BogoMIPS value is 6000.00.
\subsection{Convergence of Parareal: PINN-P vs PINO} 
\begin{figure*}[!t]
	\centering
	\includegraphics{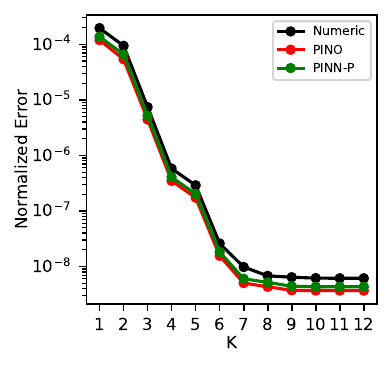}
	\hfil
        \includegraphics{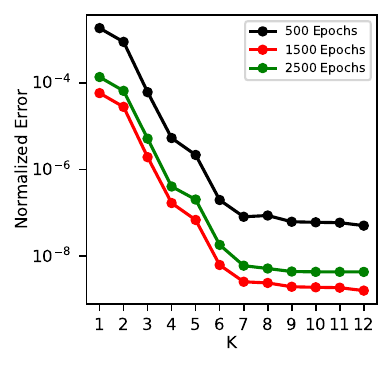}
	\caption{Convergence of Parareal against the serial fine solution with PINO, PINN-P and numerical propagator as coarse model (left). $K$ is the number of Parareal iterations. There is no significant difference in convergence between Parareal with a PINO, PINN-P or numerical coarse propagator. \rev{Convergence of PINO-Parareal against the serial fine solution if the training process is stopped after 500, 1500 or 2500 epochs (right). Longer training produces a more accurate coarse propagator which reduces error levels. However, convergence behavior is is very similar in all three cases.}}
	\label{fig:parareal_convergence_}
\end{figure*}

\rev{Figure~\ref{fig:combined_1} (left) shows the difference between the PINO-generated and the closed-form solution.}
Figure~\ref{fig:parareal_convergence_} (left) shows the normalized error of Parareal with $P_{\text{time}} = 12$ timeslices against the fine serial solution on the y-axis and the number of iterations $K$ on the x-axis.
It shows convergence for three variants of Parareal, one with a numerical coarse model, one using PINN-P as coarse model and one using PINO.
Note that because all coarse propagators run in single precision, convergence stops at a normalized error of around $10^{-8}$.
In order to provide a fair comparison of runtimes, we modified the numerical coarse propagator to also run in single precision.

There is very little difference in convergence speed between PINN-P and PINO Parareal. 
Parareal with numerical coarse models converges at the same rate but with slightly larger errors, although the differences are small.
However, the PINO uses a much smaller network than the PINN-P with fewer trainable parameters and thus is significantly faster to train, see Table~\ref{tab:pinn_vs_pino}.
The PINO is also marginally more accurate.
For comparison, the relative error of the numerical coarse propagator is $1.0 \times 10^{-2}$ while the fine propagator has a relative $l_2$-error of $2.2 \times 10^{-3}$.
Below, we report speedups for $K=1$ iterations, where Parareal will deliver a comparable discretization error as the fine propagator, as well as for $K=2$, after which Parareal has converged up to a more generic tolerance\footnote{In a use case, the discretization error of the fine propagator will not be known. Therefore, Parareal is normally run with some user-defined tolerance and $10^{-4}$ would be a reasonable value when using single precision coarse propagators.} of $10^{-4}$.

\rev{Figure~\ref{fig:parareal_convergence_} (right) shows convergence of PINO-Parareal if training of the coarse propagator is stopped at 500, 1500 or 2500 epochs.
Interestingly, the convergence behavior changes little, that is the shapes of the curves are similar, but with more training the results become more accurate.}

\begin{table}
    \caption{Size, training time and accuracy, measured as relative $l_2$-error against the closed-form solution, for the PINN-P and PINO coarse propagator.}
    \centering
    \begin{tabular}{|c|c|c|c|} \hline 
Network & Trainable parameters & Training time & Accuracy\\
PINN-P     & 1262833 & \SI{42}{\minute} & $2.3 \times 10^{-2}$ \\
PINO     & \phantom{0}832833 & \SI{10}{\minute}& $1.5 \times 10^{-2}$ \\ \hline
    \end{tabular}
    \label{tab:pinn_vs_pino}
\end{table}
    
Table~\ref{tab:components_runtime} shows the runtimes $c_{\text{fine}}$ and $c_{\text{coarse}}$ in~\eqref{eq:speedup} for the fine propagator and the numerical and PINO coarse propagator.
Upper bounds for achievable speedup according to~\eqref{eq:speedup} are $S \leq 350/113 \approx 3.1$ for Parareal-Num and $S \leq 350/2.2 \approx 159$ for Parareal-PINO.
\rev{Note that these are the runtimes per timeslice so that the total runtime for the fine serial propagator baseline would be $c_{\text{fine}} \times P_{\text{time}}$.}
\begin{table}
    \caption{Runtime of the fine propagator $\mathcal{F}$, the numerical coarse propagator and the PINO coarse propagator $\mathcal{G}$ in seconds, averaged over ten runs. Timings are done in serial execution and exclude setup times. Because maximum speedup by Parareal is bounded by $c_{\text{fine}}/c_{\text{coarse}}$, speedup for Parareal with a coarse propagator is bounded by 350/113 $\approx$ 3.1 whereas Parareal-PINO can, in theory, achieve speedup of up to 350/2.2 $\approx$ 159.1. Note that all coarse propagators run in single precision, including the numerical.}
    \centering
    \begin{tabular}{|c|c|} \hline
        \multicolumn{2}{|c|}{Runtimes \rev{per timeslice} in seconds} \\ \hline 
\rev{Fully serial fine propagator} & $350.007 \pm 0.00368$ \\
\rev{Space-parallel fine propagator ($P_{\text{space}} = 8$)} & $58.33 \pm 0.00385$  \\
\rev{Space-parallel fine propagator ($P_{\text{space}} = 16$)} & $37.82 \pm 0.002108$  \\
Numerical coarse propagator & $113.011 \pm 0.00118$ \\
PINO coarse propagator (after training) &  $\phantom{00}2.203 \pm 0.00228$ \\
PINN coarse propagator (after training) &  $\phantom{00}6.504 \pm 0.00157$ \\ \hline
    \end{tabular}
    \label{tab:components_runtime}
\end{table}

%

\subsection{Parareal-only scaling}
\begin{figure*}[!t]
    \centering
    \includegraphics{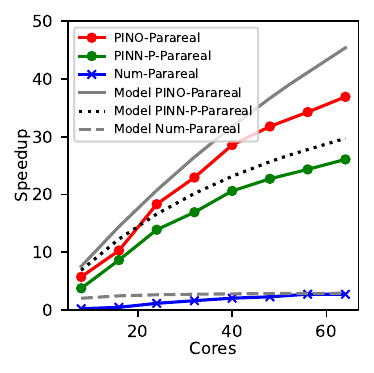}
    \hfil
    \includegraphics[scale=1]{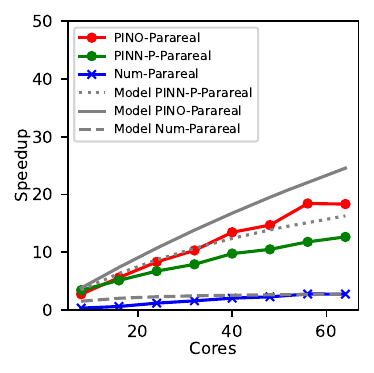}
    \caption{Parallel scaling of Parareal alone with no spatial parallelization \rev{after $K=1$ (left) and $K=2$ (right) iterations}. \rev{PINO-Parareal scale better than PINN-P-Parareal and both outperform Parareal with a numerical coarse propagator.}}
    \label{fig:parareal_scaling}
    \end{figure*}

Figure~\ref{fig:parareal_scaling} shows parallel speedup from Parareal with $K=1$ (left) or $K=2$ iterations (right), without spatial parallelization, with the number of time slices $P_{\text{time}}$ ranging from $2$ to $64$, filling the whole node.
\rev{The serially run fine propagator is used as a baseline, in line with common practice when studying performance of Parareal.}
Parareal scales reasonably well, not too far from the theoretical maximum~\eqref{eq:speedup}, even though efficiencies decrease somewhat when filling the whole node.
While both Parareal with a numerical and a PINO coarse propagator deliver speedup that is relatively close to the respective theoretical maximum, Parareal-PINO provides much better speedup throughout because of the greatly improved $c_{\text{fine}}/c_{\text{coarse}}$ ratio in~\eqref{eq:speedup}.
Note that because speedup for PINO-Parareal is limited by the second term in~\eqref{eq:speedup}, we see a much greater reduction of speedup compared to Parareal-Num when increasing $K$ from $1$ to $2$.
However, even for $K=2$, PINO-Parareal outperforms Parareal with a numerical coarse propagator and delivers more than $10$-times speedup when using $64$ timeslices.

\subsection{Space-time parallel strong scaling}
\begin{figure*}[!t]
    \centering
    \includegraphics{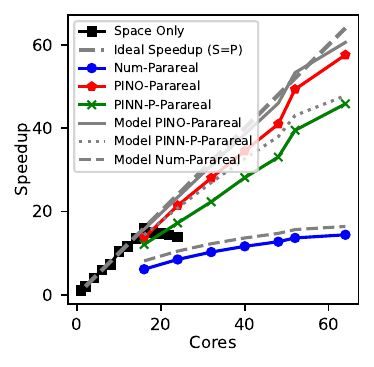}
    \hfil
    \includegraphics{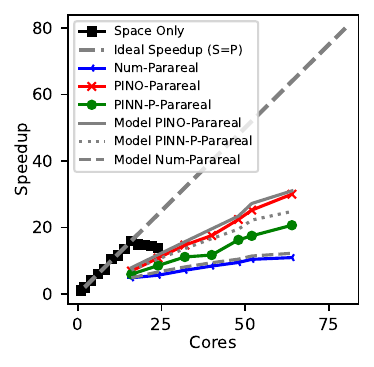}
    \caption{Space-time parallel scaling of Parareal \rev{with $K=1$ (left) and $K=2$ (right) iterations} with PINO, PINN-P and numerical coarse propagator combined with spatial parallelization using $P_{\text{space}} = 8$ cores and speedup from spatial parallelization alone. \rev{Since runtimes for the PINO coarse propagator are smaller than for the PINN-P or numerical coarse model, Parareal-PINO provides better speedup, see also the discussion after~\eqref{eq:speedup}.} }
    \label{fig:space_time_scaling}
\end{figure*}

Figure~\ref{fig:space_time_scaling} shows speedup measured against a reference simulation running the fine propagator in serial without spatial parallelization versus the total number of used cores.
Shown are speedups from a space-only parallelization using a serial fine integrator as well as space-time parallelization using Parareal with numerical and PINO coarse propagator and $K=1$ (left) and $K=2$ (right) iterations.

Spatial parallelization alone provides near-ideal speedup up to $S(P_{\text{space}}) \approx 16$ at $P_{\text{space}}=16$ cores but, after that, speedup decreases when more cores are added due to communication overheads.
Since it has been shown that it is often beneficial for space-time parallelism to use slightly fewer than the maximum number of processors up to which spatial parallelization scales~\cite{SpeckEtAl2014_Parco}, we fix $P_{\text{space}} = 8$ for all Parareal runs and increase $P_{\text{times}}$ from $2$ to $8$ for a total number of cores ranging from $P_{\text{space}} \times P_{\text{time}} = 8 \times 2 = 16$ to $8 \times 8 = 64$.    

For $K=1$, space-time Parareal with a numerical coarse propagator scales all the way up to $8 \times 8 = 64$, close to the theoretical maximum
However, total speedup, even when using the full node, is still slightly less than what space-only parallelization with $P_{\text{space}}=16$ cores can provide.
The reason is that the numerical coarse propagator is too expensive compared to the fine.
In such a case, Parareal will not provide any benefits.
By contrast, Parareal-PINO provides close to ideal speedup all the way up to the full node.
The combination of Parareal-PINO with spatial parallelization extends scaling from a total speedup of around $16$ to a total speedup of around $50$ by being able to efficiently utilize the full $64$ cores of the node.

For $K=2$ iterations, the speedup of Parareal-PINO in Figure~\ref{fig:space_time_scaling} is noticeably reduced due to the second bound in~\eqref{eq:speedup}.
Because speedup from Parareal is bounded by $P_{\text{time}} / K$ and thus parallel efficiency by $1 / K$, we now also see speedup that is farther away from ideal speedup.
Here, around $40$ cores are required for space-time parallelization to become faster than purely spatial parallelization.
However, when using the full node, space-time parallelization using Parareal-PINO still extends scaling, providing speedup of around $30$ in contrast to $16$ from spatial parallelization alone.

\begin{figure*}[!t]
    \centering
    \includegraphics{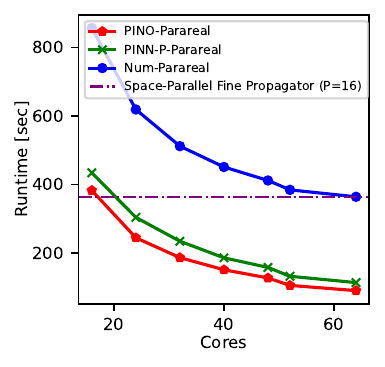}
    \hfil
    \includegraphics{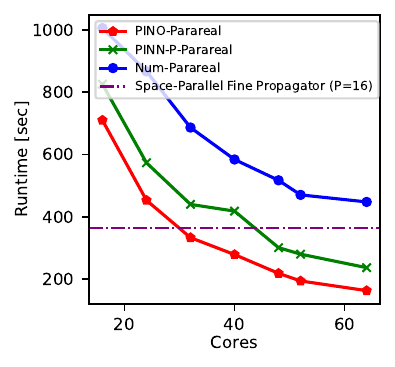}
   \caption{\rev{Runtimes of Parareal with numerical, PINO and PINN-P coarse propagator for $K=1$ (left) and $K=2$ (right) iterations. The runtime of the space-parallel time-serial fine propagator using $P_{\text{space}}=16$ processes is shown as horizontal lines.}
}
    \label{fig:parareal_runtime_comparison}
\end{figure*}

\rev{Figure~\ref{fig:parareal_runtime_comparison} shows runtimes of Parareal with PINO, PINN-P and numerical coarse propagator and $P_{\text{space}}=8$ processes for parallelization in space.
For comparison, the serial fine propagator with no parallelization in space requires $\approx$ 3150.63 seconds while at peak spatial parallelization using $P_{\text{space}}=16$ processes takes $\approx 364$ seconds.
The latter value is indicated in both plots as a dashed-dotted horizontal line.
In both cases, Parareal with a numerical coarse propagator struggles to reduce runtimes below what space parallelization alone can provide.
However, both PINN-P-Parareal and PINO-Parareal, in combination with spatial parallelization, can reduce runtimes further.
In both cases, runtimes are smaller for PINO-Parareal than for PINN-P-Parareal.}
\subsection{Generalization to different resolution}
\begin{figure}[!t]
    \centering
    \includegraphics{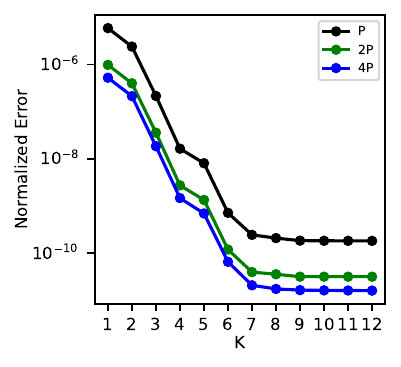}
    \caption{Convergence of PINO-Parareal when increasing spatial and temporal resolution and thus problem size as well as the number of timeslices together (``weak scaling''). The PINO coarse propagator generalizes well to different spatial and temporal resolutions.}
    \label{fig:parareal_scaling_convergence}    
\end{figure}

Figure~\ref{fig:parareal_scaling_convergence} shows convergence of PINO-Parareal against the serial fine solution in a weak scaling type test.
We start with $\Delta x = \Delta y = 1 $ spatial resolution and $\Delta \tau = 0.1$ temporal resolution in the coarse and $\Delta \tau/3$ in the fine propagator using $P_{\text{time}} = P = 2$ time slices.
We twice double both spatial and temporal resolution and the number of timeslices to $P_{\text{time}} = 2P = 4$ and $P_{\text{time}} = 4P = 8$.
The PINO remains the same as before and is not retrained — \rev{this is possible because PINOs are inherently mesh-invariant, as they learn mappings between function 
spaces rather than finite-dimensional vectors. 
By contrast, traditional neural networks operate on fixed-size inputs and outputs.}

Convergence of PINO-Parareal remains fast and error levels even decrease slightly as resolution increases.
For a fixed error tolerance, the number of required Parareal iterations therefore remains constant as resolution and problem size increase.
This means that, \rev{for the Black-Scholes equations}, PINO-Parareal can deliver good weak scaling if implemented efficiently because its computational cost will not grow with problem size.
\rev{While similar behavior can be expected for other parabolic PDES where Parareal typically performs well, this will not be true for 
hyperbolic-style problems, where numerical diffusion becomes weaker as $\Delta x \to 0$, which will make Parareal converge worse~\cite{GanderVandewalle2007_SISC}.}

\subsection{Generalization \rev{of Parareal-PINO} to different model parameters}
\begin{figure*}[!t]
    \centering
    \includegraphics{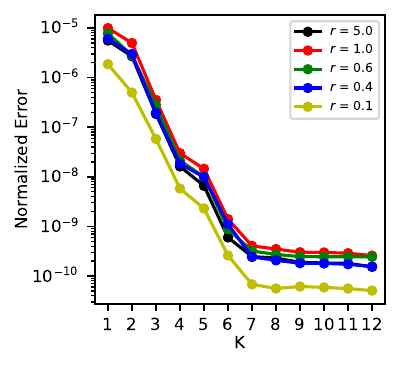}
    \includegraphics{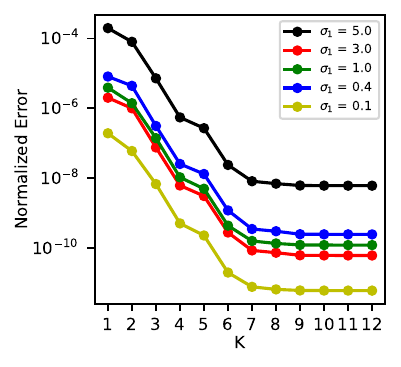}\\
    \includegraphics{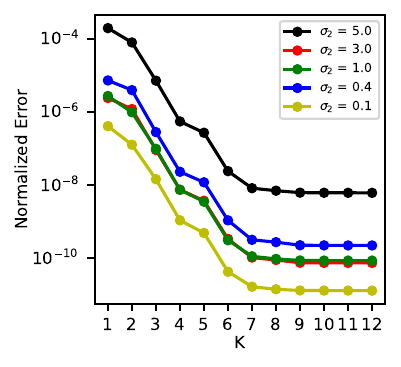}
    \includegraphics{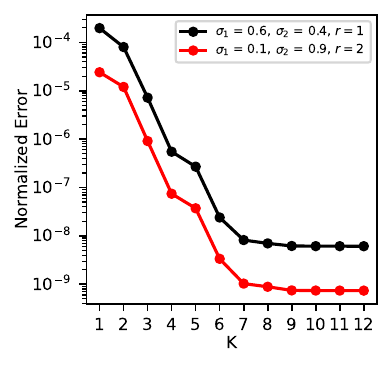}
    \caption{Convergence of PINO-Parareal for model parameters different to what the PINO was trained for. The values used for training are $\sigma_1 = 0.2$, $\sigma_2 = 0.3$ and $r = 1$. The model is very robust against changes in $r$. Changes in volatilities $\sigma_1$ and $\sigma_2$ have a more discernible impact but, even when far exceeding the training values, PINO-Parareal still converges fast, even though error levels are higher.}
        \label{fig:parareal_convergence_diff}  
\end{figure*}

The PINO was trained for model parameters $\sigma_1 = 0.3$, $\sigma_2 = 0.3$ and $r = 1$ and resolution $\Delta x = \Delta y = 1$ and $\Delta \tau = 0.1$.
To assess how well PINO-Parareal generalizes and how robust its convergence is against changes in model parameters, we show the convergence error of PINO-Parareal against the fine serial solution for a range of different values for $r$, $\sigma_1$ and $\sigma_2$ in Figure~\ref{fig:parareal_convergence_diff}.
The PINO is not retrained or otherwise modified - \rev{therefore, PINO-Parareal uses a coarse propagator that was trained for different parameters than the problem it is applied to.}
The rate of convergence remains mostly invariant but error levels change as the parameters change.
Convergence is very robust against variations in $r$.
Even for $r=5$ the difference in convergence compared to the training value of $r=1$ is small.
Sensitivities are larger for changes in the volatilities $\sigma_1$ and $\sigma_2$.
Interestingly, values of $\sigma_1, \sigma_2$ that are slightly smaller than the trained value improve convergence, possibly because the dynamics of the problem become less rapid.
However, in the limit $\sigma_1, \sigma_2 \to 0$, only the first order derivatives are left in~\eqref{eq:bse_change_variable}.
For very small volatilities, the equation degenerates into an advection-type problem and we see the well-documented non-monotonic convergence of Parareal for transport problems~\cite{GanderVandewalle2007_SISC}.
These results are not shown here but can be generated with the provided code.
For values of $\sigma_1, \sigma_2$ larger than the training value, error levels increase until, at $\sigma_1 = 5$ or $\sigma_2=5$, they are about one order of magnitude larger.
However, at this point, the training values have been exceeded by a factor of $5/0.2 = 25$ for $\sigma_1$ and $5/0.3 \approx 17$ for $\sigma_2$, so that some degeneration of performance is to be expected.

For our test cases, we used a range of values for both volatility $\sigma$ and the risk-free interest rate $r$: 0.1, 0.4, 1, 3, and 5. 
While $sigma$ and $r$ values of 5 are large compared to typical market conditions — where volatilities usually range from 15\% to 60\% and interest rates are often below 2\% — 
these values were intentionally chosen to test the models robustness and evaluate how well the PINO generalizes across a wide parameter space, including edge cases. 
This approach helps ensure the model can handle both realistic and high-stress scenarios.

\section{Summary}\label{sec:discussion}
We investigate convergence and parallel scaling of the Parareal parallel-in-time method using a Fourier Neural Operator as a coarse model.
As a benchmark problem, we use the two-asset Black-Scholes equation from computational finance.
The PINO provides accuracy comparable to a numerical coarse model and a previously studied coarse model based on a physics-informed neural network propagator (PINN-P).
However, the PINO takes only about a quarter of the training time of the PINN-P, making it a better choice.
Because evaluating the PINO, once trained, is about a factor of fifty faster than running the numerical coarse model \rev{and about three times faster than the PINN-P}, Parareal-PINO greatly relaxes the bound on speedup for Parareal given by the ratio of runtimes of the fine to the coarse propagator.
Whereas the speedup of standard Parareal is bounded by 3.1, the much faster PINO coarse model theoretically allows for Parareal speedups up to 159.
We perform scaling tests of Parareal alone and of a combination of Parareal with spatial parallelization.
In both cases, Parareal-PINO significantly outperforms standard Parareal.
Furthermore, we demonstrate that a combined space-time parallelization using Parareal-PINO scales beyond the saturation point of spatial parallelization alone.
Whereas the latter saturates at a speedup of around $16$ on as many cores, the former scales to the full node with $64$ cores, providing a total speedup of almost $60$ after one and $30$ after two Parareal iterations.  
We also show that convergence of PINO-Parareal is robust against changes in parameters or resolution without retraining the model.

\section*{Data availability.}
The code for this paper is available from Zenodo~\cite{IbrahimRepo}.

\begin{acks}
  This project has received funding from the European High-Performance Computing Joint Undertaking (JU) under grant agreement No 101118139. The JU receives support from the European Union's Horizon Europe Programme.
\end{acks}
\bibliographystyle{ACM-Reference-Format}
\bibliography{pint,refs}
\end{document}